\documentclass[12pt]{amsart}
\usepackage{amscd,amssymb,float}
\newtheorem{theorem}{Theorem}[section]
\newtheorem{lemma}[theorem]{Lemma}

\newtheorem{proposition}[theorem]{Proposition}

\theoremstyle{definition}

\theoremstyle{remark}

\theoremstyle{conj}

\numberwithin{equation}{section}

\begin{document}
\title{Bounding sectional curvature along a K\"ahler-Ricci flow  }
\author[W.Ruan]{Wei-Dong Ruan}
\address{Department of Mathematical Sciences, Korea
Advanced Institute of Science and Technology, Daejeon, Republic of
Korea}
  \email{ruan@math.kaist.ac.kr}
\author[Y. Zhang]{Yuguang Zhang}
\address{Department of Mathematics, Capital Normal University,
Beijing, P.R.China  }\address{Department of Mathematical Sciences,
Korea Advanced Institute of Science and Technology, Daejeon,
Republic of  Korea}
 \email{zhangyuguang76@yahoo.com}
\author[Z. Zhang]{Zhenlei Zhang}
\address{Nankai Institute of Mathematics,
Weijin Road 94, Tianjin 300071, P.R.China}
 \email{zhleigo@yahoo.com.cn}

\begin{abstract}
If a normalized K\"{a}hler-Ricci flow $g(t),t\in[0,\infty),$ on a
compact K\"{a}hler manifold $M$, $\dim_{\mathbb{C}}M=n\geq 3$,  with
positive first Chern class satisfies $g(t)\in 2\pi c_{1}(M)$ and has
 curvature operator uniformly bounded in  $L^{n}$-norm,   the curvature operator
will also be  uniformly bounded along the flow. Consequently the
flow will converge along a subsequence to a K\"{a}hler-Ricci
soliton.

\end{abstract}
\maketitle

\section{Introduction}
On a compact K\"ahler manifold $M$, $\dim_{\mathbb{C}}M=n$, with the
first Chern class $c_{1}(M)>0$, the normalized K\"ahler-Ricci flow
equation  is
\begin{equation}\label{1.1}\partial_{t}g(t)=-Ric(g(t))+g(t),
\end{equation} for a family of K\"ahler metrics $g(t)\in 2\pi c_1(M)$, where, for brevity, $g(t)$
 denotes  either  K\"ahler
 metrics or  K\"ahler forms depending on the context.   In
\cite{Cao}, it is proved that  a solution $g(t)$ of (\ref{1.1})
exists
 for all  $t\in [0, \infty)$. Perelman (cf. \cite{ST})) has proved some
 important properties for the solution $g(t)$, $t\in [0, \infty)$,  of
 (\ref{1.1}): there exist constants $C>0$ and  $\kappa >0$ independent of $t$ such
 that
 \begin{itemize}\label{00}
  \item[(1)] $|R(g(t))|< C$, and  $\text{diam}_{g(t)}(M)<C$,
  \item[(2)] $(M, g(t))$ is $\kappa$-noncollapsed, i.e. for any
$r< 1$,  if $|R(g(t))|\leq r^{-2}$ on a metric ball $B_{g(t)}(x,r)$,
then
\begin{equation}\label{2.5}\text{Vol}_{g(t)}(B_{g(t)}(x,r))\geq \kappa r^{2n}.
\end{equation}

 \end{itemize}
 In a recent preprint \cite{S3},
 Sesum has proved that, if $n\geq
 3$, assuming  the Ricci curvatures  $|Ric(g(t))|<C$ and the integral of   curvature
 operators $
\int_{M}|Rm(g(t))|^{n}dv_{t}\leq C,$ for a constant $C$ in-dependent
of $t$, then the curvature
 operators are uniformly bounded.
   In this note , we will  show  that the hypothesis of
bounded Ricci curvature can be removed.

\vskip 5mm
\begin{theorem} Let $g(t)$, $t\in [0, +\infty)$, be a solution of the
normalized K\"ahler-Ricci flow (\ref{1.1}) on a compact K\"ahler
 manifold $M$ with $c_1(M)>0$ and initial metric  $g (0)\in 2\pi
c_1(M)$. Assume that $dim_{\mathbb{C}}M=n\geq 3$.
 If the $L^{n}$-norms of curvature operators  are  uniformly bounded by a constant
 $C$, i.e.
$$
\int_{M}|Rm(g(t))|^{n}dv_{t}\leq C,$$  then there  exists a constant
$0< \bar{C}<\infty$ such that
\[
\sup_{M\times[0,\infty)}|Rm(g(t))|\leq\bar{C}.
\]
Consequently the flow will converge along a subsequence to a
K\"{a}hler-Ricci soliton.
\end{theorem}
\vskip 5mm

From this theorem,   it is a direct  consequence of Hamilton's
compactness theorem (c.f.
 \cite{H1}) that, for any $t_{k}\rightarrow \infty$, a
subsequence of $( M, g(t_{k}+t))$, $t\in [0, 1]$,  converges
smoothly to  $(X, h(t))$, $t\in [0, 1]$, where $X$ is  a  compact
 complex  manifold, and   $\{h(t)\}$, $t\in [0, 1]$,  is a family of K\"ahler metrics that    satisfies
 the K\"ahler-Ricci flow  equation. Furthermore, from the arguments in  the proof of Theorem 12 in \cite{ST}, $h(t)$,
$t\in [0,1]$, satisfies the K\"ahler-Ricci soliton equation, i.e.
there is a holomorphic vector field  $v$ on $X$ such that
$$Ric(h)-h=\mathcal{L}_{v}h.$$

In \cite{S3}, Sesum conjectured that Theorem 1.1, as stated for
$n\geq 3$,  is also true for $n=2$. By the classification theory of
complex surface, the only compact K\"ahler surfaces with $c_1(M)>0$
are diffeomorphic to $\mathbb{CP}^{2}\sharp
l\overline{\mathbb{CP}}^{2}$,  $0\leq l\leq8$, and
$\mathbb{CP}^{1}\times \mathbb{CP}^{1}$. By \cite{T}, each of
$\mathbb{CP}^{1}\times \mathbb{CP}^{1}$ and $\mathbb{CP}^{2}\sharp
l\overline{\mathbb{CP}}^{2}$,  $3\leq l\leq8$ or $l=0$, admits a
K\"ahler-Einstein metric. In \cite{C2} and \cite{Ko}, it is shown
that $\mathbb{CP}^{2}\sharp \overline{\mathbb{CP}}^{2}$ admits a
non-trivial  K\"ahler-Ricci soliton  metric. Wang and Zhu \cite{WZ}
showed the same result later for $\mathbb{CP}^{2}\sharp
2\overline{\mathbb{CP}}^{2}$. By \cite{TZ}, on a compact  K\"ahler
surface $M$ with $c_1(M)>0$, if the initial metric $g(0)$ is
invariant under a one-parameter group obtained from a K\"ahler-Ricci
soliton  metric on $M$, the curvatures stay uniformly bounded along
the flow. The only remaining case  is when  $M$  is a complex
surface diffeomorphic to $\mathbb{CP}^{2}\sharp
\overline{\mathbb{CP}}^{2}$ or $\mathbb{CP}^{2}\sharp
2\overline{\mathbb{CP}}^{2}$, with an initial metric $g(0)$ without
any symmetry. In \cite{FZ}, it is
 proved  that the K\"ahler-Ricci flow on $\mathbb{CP}^{2}\sharp
\overline{\mathbb{CP}}^{2}$ or $\mathbb{CP}^{2}\sharp
2\overline{\mathbb{CP}}^{2}$ converges to an  orbifold  in the
Gromov-Hausdorff sense.

 By using the
method in the proof of Theorem 1.1,  we can give a  different proof
of the convergence of the  K\"ahler-Ricci flow on $\mathbb{CP}^{2}$
(Theorem 3.3),
 which is already implied by \cite{TZ}.   In a very recent
preprint \cite{CW}, Chen and Wang claimed  that the bounding  of
curvatures along the K\"ahler-Ricci  flow on a toric Fano surface
$M$ (including $\mathbb{CP}^{2}\sharp \overline{\mathbb{CP}}^{2}$
and $\mathbb{CP}^{2}\sharp 2\overline{\mathbb{CP}}^{2}$)  could be
proved by using the fact that $M$ is a  toric manifold.

There is also an analogy to Theorem 1.1 in the real Ricci flow case.

\begin{theorem}
Let $g(t),t\in[0,T),$ be a solution to the Ricci flow, normalized or
not, on a closed odd dimensional manifold $M$ with $T<\infty$.
Suppose that
$$\int_{M}|Rm(g(t))|^{n/2}dv_{t}\leq C, \ \ \ \ {\rm and } \ \ \ \ \ |R(g(t))|\leq C$$ for a
constant $C<\infty$ independent of $t$, where
$n=dim_{\mathbb{R}}(M)$. Then there is another constant
$\bar{C}<\infty$ such that
$$\sup_{M\times[0,T)}|Rm(g(t))|\leq\bar{C},$$
and so the flow can be extended over $T$.
\end{theorem}

The organization of the paper is as follows: In $\S$2,  we prove
Theorem 1.1. In $\S$3,  we give some remarks for K\"ahler-Ricci on
Fano surfaces.  Then we prove Theorem 1.2 in $\S$4.

\vskip 10mm

\noindent {\bf Acknowledgement:} The authors would like to  thank Bing
Wang for informing them the paper \cite{CW}.

\section{Proof of Theorem 1.1}

Let $g(t)$, $t\in [0, +\infty)$, be a solution of the normalized
K\"ahler-Ricci flow (\ref{1.1}) on a compact K\"ahler manifold $M$,
  $\dim_{\mathbb{C}}M=n$,  with $c_1(M)>0$ and initial metric  $g
(0)\in 2\pi c_1(M)$.
 Assume that $dim_{\mathbb{C}}M=n\geq 3$, and that
$$
\int_{M}|Rm(g(t))|^{n}dv_{t}\leq C,$$ for  a constant $C$
independent of $t$.   Perelman (cf. \cite{ST})) has proved that
there exist constants $C>0$, $\kappa
>0$ independent of $t$ such
 that
\begin{itemize}\label{00}
  \item[(1)] $|R(g(t))|< C$, and  $\text{diam}_{g(t)}(M)<C$,
  \item[(2)] $(M, g(t))$ is $\kappa$-noncollapsed, i.e. for any
$r< 1$,  if $|R(g(t))|\leq r^{-2}$ on a metric ball $B_{g(t)}(x,r)$,
then
\begin{equation}\label{2.1}\text{Vol}_{g(t)}(B_{g(t)}(x,r))\geq \kappa r^{2n}.
\end{equation}

  \end{itemize}

    The proof of  Theorem 1.1  relies on the following  theorem due to Gang Tian:
\begin{theorem}[Theorem 2 in \cite{T2}] Let $(N, J, g)$ be a
complete non-compact Ricci-flat K\"ahler manifold with
$dim_{\mathbb{C}}M=n\geq 3$,
$$\int_{N}|Rm(g)|^{n}dv_{g}< C<\infty, \ \ \ \ \ \ \ \rm and $$
 $$\text{Vol}_{g}(B_{g}(x, r))\geq \kappa r^{2n},$$ for any $r>0$, where $C$
and $\kappa$ are constants. Then $M$ is a resolution of
$\mathbb{C}^{n}/\Gamma$ where $\Gamma$ is a
  finite group $\Gamma \subset SU(n)$, which  acts  on
$\mathbb{C}^{n}\backslash \{0\}$ freely, i.e. there is a holomorphic
map $\pi: N\longrightarrow \mathbb{C}^{n}/\Gamma$ such that $\pi:
N\backslash \pi^{-1}(0)\longrightarrow \mathbb{C}^{n}\backslash
\{0\}/\Gamma$ is bi-holomorphic.
\end{theorem}

 The assumptions of  Euclidean volume growth and
$dim_{\mathbb{C}}M=n\geq 3$ not mentioned explicitly in Theorem 2 in
\cite{T2} seem to be necessary.
 Let's recall several main steps in the proof of
Theorem 2.1. First, an estimate for the decrease  of the sectional
curvature of $g$ is obtained by assuming the Euclidean volume
growth,  bounded  $L^{n}$-norm of curvature operator, and the
Ricci-flat metric (See Lemma 4.1 in \cite{T2}). Then, for proving
Theorem 2 of \cite{T2}, one  needs  Lemma 3.4 and Lemma 3.3 of
\cite{T2}, which  have  the hypothesis $dim_{\mathbb{C}}M=n\geq 3$.
The main tool there was Kohn's estimate for $\square_{b}$-operators
that works only for $n\geq 3$ (See \cite{T2} for details).

\begin{proof}[Proof of Theorem 1.1]

Suppose otherwise, there exists  a sequence of times
$t_{k}\rightarrow\infty$, and a sequence of  points $x_{k}\in M$
such that
$$Q_{k}=|Rm(g(t_{k}))|(x_{k})=\sup_{M\times[0,t_{k}]}|Rm(g(t))|\rightarrow\infty.$$
Consider the sequence $\widetilde{g}_{k}(t),t\in[-Q_{k}t_{k},0]$,
where  $\widetilde{g}_{k}(t)=Q_{k}g(Q_{k}^{-1}t+t_{k}))$ and satisfy
\begin{equation}\label{2.2}
\partial_{t}\widetilde{g}_{k}(t)=-Ric(\widetilde{g}_{k}(t))+Q_{k}^{-1}\widetilde{g}_{k}(t), \end{equation}
\begin{equation}\label{2.3}\sup_{M\times[-Q_{k}t_{k},0]}|Rm(\widetilde{g}_{k}(t))|\leq 1, \ \
\ \ \ {\rm and}  \ \ \ \ \ |Rm(\widetilde{g}_{k}(0))|(x_{k})=
1.\end{equation} By Perelman's estimate, we obtain that
 \begin{equation}\label{2.4}|R(\widetilde{g}_{k}(t))|< CQ_{k}^{-1}\longrightarrow 0,\end{equation}
when $k\rightarrow \infty$, and,
 for any
$r< CQ_{k}^{\frac{1}{2}}$,  $x\in M$,
 \begin{equation}\label{2.5}\text{Vol}_{\widetilde{g}_{k}(t)}(B_{\widetilde{g}_{k}(t)}(x,r))\geq
\kappa r^{2n}. \end{equation}

 By Hamilton's
compactness theorem (c.f. Appendix E in \cite{KL}), by passing to a
subsequence, $\{(M, J, \widetilde{g}_{k}(t),  x_{k})\}$, $t\in
[-1,0]$,  converges smoothly  to a family of pointed  complete
Riemannian manifold $(N, J_{\infty}, g_{\infty}(t), x_{\infty})$,
$t\in [-1,0]$, where $g_{\infty}(t)$, $t\in [-1,0]$, satisfies the
Ricci-flow equation $\partial_{t}g_{\infty}(t)=-Ric(g_{\infty}(t))$.
Particularly, for any $k\gg 1$ and $r>0$, there is an embedding
$F_{k,r}: B_{g_{\infty}}(x_{\infty},r)\longrightarrow M$ such that
$F_{k,r}^{*}\widetilde{g}_{k}(0)$ converges smoothly to
$g_{\infty}$, and $dF_{k,r}^{-1}JdF_{k,r}$ converges smoothly to an
almost complex structure $J_{\infty}$, where
$g_{\infty}=g_{\infty}(0)$.  Actually, $J_{\infty}$ is integrable,
and $g_{\infty}$ is a K\"ahler metric of $J_{\infty}$ (c.f.
\cite{Ru}). From (\ref{2.3}) and (\ref{2.5}), we obtain that
\begin{equation}\label{2.6}|Rm(g_{\infty})|\leq|Rm(g_{\infty})|(x_{\infty})=
1, \end{equation}
 and, for any
$r>0$ and $x\in N$,
\begin{equation}\label{2.7}\text{Vol}_{g_{\infty}(t)}(B_{g_{\infty}(t)}(x,r))\geq \kappa
r^{2n}. \end{equation} By  (\ref{2.4}),   $R(g_{\infty}(t))\equiv
0$, which  implies that $Ric(g_{\infty}(t))\equiv 0$ since
$g_{\infty}(t)$ is a solution of the Ricci-flow equation.
 From  the smooth  convergence,
$$\int_{N}|Rm(g_{\infty})|^{n}dv_{g_{\infty}}\leq \limsup_{k\longrightarrow \infty}
\int_{M}|Rm(g(t_{k}))|^{n}dv_{g(t_{k})}\leq C<\infty.$$

Thus  $(N, J_{\infty}, g_{\infty})$ is a complete Ricci-flat
K\"ahler manifold with  Euclidean volume growth, and  $L^{n}$-norm
of curvature operator bounded.  By (2.6),    $
  g_{\infty}$
is not a   flat metric. Note  that $dim_{\mathbb{C}}M=n\geq 3$. By
Theorem 2.1, $N$ is a resolution of $\mathbb{C}^{n}/\Gamma$ where
$\Gamma$ is a
  finite group $\Gamma \subset SU(n)$, which  acts  on
$\mathbb{C}^{n}\backslash \{0\}$ freely, i.e. there is a holomorphic
map $\pi: N\longrightarrow \mathbb{C}^{n}/\Gamma$ such that $\pi:
N\backslash \pi^{-1}(0)\longrightarrow \mathbb{C}^{n}\backslash
\{0\}/\Gamma$ is bi-holomorphic. If $\Gamma$ is trivial, i.e.
$\Gamma=\{e\}$, then $(N,g_{\infty})$ is isometric to
$\mathbb{R}^{2n}$ by Theorem 3.5 in \cite{An1}, which contradicts
(2.1.5). Thus $\Gamma$ is non-trivial,  and $V=\pi^{-1}(0)$ is a
compact analytic subvariety of $(N, J_{\infty})$ with
$0<dim_{\mathbb{C}}V=m<n$ (See \cite{GH} for the definition of
$dim_{\mathbb{C}}V$). We obtain that
$$\int_{V}g_{\infty}^{m}=m!Vol_{g_{\infty}}(V)> 0.$$

Note that, for any $k$, $F_{k,r}(V)$ is  a cycle, and  defines a
homology class $[F_{k,r}(V)]\in H_{2m}(M, \mathbb{Z})$.   By the
smooth  convergence of $F_{k,r}^{*}\widetilde{g}_{k}$, for any
$\varepsilon
>0$, there is a $k_{0}> 0$ such that, for any $k\geq k_{0}$,
\[
\left|\int_{V}g_{\infty}^{m}-\int_{F_{k,r}(V)}\widetilde{g}_{k}^{m}\right|\leq
\varepsilon,
\]
where $\widetilde{g}_{k}=\widetilde{g}_{k}(0)$, and $r\gg 1$ such
that $V\subset B_{g_{\infty}}(x_{\infty}, r)$.
 As $g_{k}=g(t_{k})\in 2\pi c_1(M)$,
$$\int_{F_{k,r}(V)}\widetilde{g}_{k}^{m}
=Q_{k}^{m}\int_{F_{k,r}(V)}g_{k}^{m}=Q_{k}^{m}(2\pi)^{m}\int_{F_{k,r}(V)}c_1^m(M).$$
  By taking
$\varepsilon=\frac{1}{2}\int_{V}g_{\infty}^{m}$ and $k\gg k_{0}$, we
obtain that
$$0<\frac{1}{2}Q_{k}^{-m}(2\pi)^{-m}\int_{V}g_{\infty}^{m}\leq
\int_{F_{k,r}(V)}c_1^m(M)\leq
\frac{3}{2}Q_{k}^{-m}(2\pi)^{-m}\int_{V}g_{\infty}^{m} <1.$$ Since
$0\neq c_1^m(M)\in H^{2m}(M,\mathbb{Z})$, and $[F_{k,r}(V)]\in
H_{2m}(M,\mathbb{Z})$, we have
$$\int_{F_{k,r}(V)}c_1^m(M)\in \mathbb{Z}.$$ It is a
contradiction. We obtain  that
\[
\sup_{M\times[0,\infty)}|Rm(g(t))|\leq\bar{C},
\] for a constant $\bar{C}>0$.

Now, by  Hamilton's compactness theorem (c.f.
 \cite{H1}), for any $t_{k}\rightarrow \infty$, a
subsequence of $( M, g(t_{k}+t))$, $t\in [0, 1]$,  converges
smoothly to a family of  compact K\"ahler manifolds  $(X, h(t))$,
$t\in [0, 1]$, where $h(t)$   satisfies
 the K\"ahler-Ricci flow  equation.  Actually,   $h(t)$,
$t\in [0,1]$, satisfies the K\"ahler-Ricci soliton equation from the
arguments in  the proof of Theorem 12 in \cite{ST}.

\end{proof}

\section{Remarks for K\"ahler
 surfaces}
 Let   $g(t)$, $t\in [0, +\infty)$, be   a solution of the normalized
K\"ahler-Ricci flow (\ref{1.1}) on a compact K\"ahler surface  $M$,
i.e. $\dim_{\mathbb{C}}M=2$,  with $c_1(M)>0$ and initial metric  $g
(0)\in 2\pi c_1(M)$.

 \begin{lemma}
 The $L^{2}$-norms of curvature operators of $g(t)$
are bounded along the flow, i.e. there is a constant $C>0$
independent of $t$ such that
\begin{equation}\label{2.8} \int_{M}|Rm(g(t))|^{2}dv_{t}\leq C.\end{equation}
 \end{lemma}

 \begin{proof} Since, for any $t\in [0,
\infty)$, $(M, g(t))$ is  a K\"ahler
 surface, we have $\int_{M}c_{1}^{2}(M)=2\chi(M)+3\tau (M)$,
 $R^{2}(g(t))=24|W^{+}(g(t))|^{2}$, and Gauss-Bonnet-Chern formula
 and Hirzebruch formula
 $$\chi(M)=\frac{1}{8\pi^{2}}\int_{M}(\frac{R^{2}}{24}+|W^{+}|^{2}+|W^{-}|^{2}-\frac{1}{2}|Ric\textordmasculine|^{2})dv_{t},
$$ $$\tau(M)=\frac{1}{12\pi^{2}}\int_{M}(|W^{+}|^{2}-|W^{-}|^{2})dv_{t}
$$ (c.f.  \cite{Be}),   where $Ric\textordmasculine
 =Ric(g(t))-\frac{R}{4}g(t)$,  $W^{\pm}(g(t))$ are the self-dual and
anti-self-dual Weyl tensors  of $g(t)$, and $\chi(M)$ (respectively
$ \tau (M)$) is the Euler number (respectively signature) of $M$.
Then we obtain that
$$\int_{M}|Ric\textordmasculine|^{2}dv_{t}=\int_{M}\frac{R^{2}}{4}dv_{t}-8\pi^{2}c_{1}^{2}(M),
\ \  {\rm and } $$
\begin{equation} \ \
\int_{M}|W^{-}|^{2}dv_{t}=\int_{M}\frac{R^{2}}{24}dv_{t}-12\pi^{2}\tau(M).
\end{equation} Note that  \[ Rm
(g(t))=\begin{pmatrix}W^{+}+\frac{R}{12} & Ric\textordmasculine
\\ Ric\textordmasculine  & W^{-}+\frac{R}{12}\end{pmatrix}. \] Thus
 we obtain  that \begin{eqnarray*}
\int_{M}|Rm(g(t))|^{2}dv_{t}&= &
\int_{M}(\frac{R^{2}}{24}+|W^{+}|^{2}+|W^{-}|^{2}+2|Ric\textordmasculine|^{2})dv_{t}
\\ & = & \int_{M}\frac{5R^{2}}{8}dv_{t}-16\pi^{2}c_{1}^{2}(M)
-12\pi^{2}\tau(M) .\end{eqnarray*}   Hence, by Perelman's estimate
for scalar curvatures, we obtain (3.1). \end{proof}

  Unfortunately,  our arguments in the
proof of Theorem 1.1 can not be generalized to this  case, even for
$\mathbb{CP}^{2}\sharp
 \overline{\mathbb{CP}}^{2}$.   The essential point in  the
 proof of Theorem 1.1 is that, in any Asymptotically Locally Euclidean Ricci-flat
 K\"ahler manifold $N$ of $\dim_{\mathbb{C}}N\geq 3$,  we can find a    non-trivial  class  $[V]\neq 0\in
H_{2m}(N,\mathbb{Z})$ with $m \geq 1$. However,
 there are ALE Ricci-flat
 K\"ahler surfaces without  such homology classes.   For
 example, there is an ALE  Ricci-flat
 K\"ahler metric $h$ on $T^{*}\mathbb{RP}^{2}$ whose Betti numbers satisfy $b_{2}=b_{3}
 =b_{4}=0$. Actually, the
 universal covering space of $T^{*}\mathbb{RP}^{2}$ with the pull-back metric is the  Eguchi-Hanson space (c.f. \cite{D}), which is diffeomorphic
 to $T^{*}S^{2}$.

  \begin{proposition} Assume that  $M$ is   diffeomorphic to $\mathbb{CP}^{2}\sharp
\overline{\mathbb{CP}}^{2}$.  If there is a  sequence of times
$t_{k}\rightarrow\infty$ such that
$$Q_{k}=|Rm(g(t_{k}))|(x_{k})=\sup_{M\times[0,t_{k}]}|Rm(g(t))|\rightarrow\infty,$$
where  $x_{k}\in M$,  then  a subsequence of $(M, Q_{k}g(t_{k}),
x_{k})$ converges smoothly to an ALE Ricci-flat
 K\"ahler surface $(N, g_{\infty}, x_{\infty})$ in the pointed
 Gromov-Hausdorff sense. Furthermore,
 the fundamental group $ \pi_{1}(N)$ of $N$ is a  non-trivial finite   group.
 \end{proposition}

\begin{proof}
Let $\widetilde{g}_{k}=Q_{k}g(t_{k})$.  By  the same arguments as in
the proof of Theorem 1.1, by passing to a subsequence, $\{(M, J,
\widetilde{g}_{k}, x_{k})\}$ converges smoothly to a  complete
Ricci-flat K\"ahler surface   $(N, J_{\infty}, g_{\infty},
x_{\infty})$, i.e.  for any $k\gg 1$ and $r>0$, there is an
embedding $F_{k,r}: B_{g_{\infty}}(x_{\infty},r)\longrightarrow M$
such that $F_{k,r}^{*}\widetilde{g}_{k}$ converges smoothly to
$g_{\infty}$, and $dF_{k,r}^{-1}JdF_{k,r}$ converges smoothly to
$J_{\infty}$. Furthermore,  $(N, J_{\infty}, g_{\infty})$  satisfies
that,
 for any
$r>0$ and $x\in N$,
$$\text{Vol}_{g_{\infty}}(B_{g_{\infty}}(x,r))\geq \kappa
r^{4},
$$
$$\int_{N}|Rm(g_{\infty})|^{2}dv_{g_{\infty}}\leq C<\infty,$$ $$  {\rm and} \ \
 \ \ \  |Rm(g_{\infty})|\leq|Rm(g_{\infty})|(x_{\infty})= 1. $$
 Thus, by Theorem 1.5 in \cite{BKN}, $(N, J_{\infty}, g_{\infty})$
is an Asymptotically Locally Euclidean Ricci-flat K\"ahler surface.
 Since $ g_{\infty}$ is not flat, it is easy to see that
 the fundamental group $ \pi_{1}(N)$ of $N$ is a finite group (c.f.
\cite{An1}).

If  $ \pi_{1}(N)= \{1\}$, $(N, J_{\infty}, g_{\infty})$ is an ALE
hyper-K\"ahler 4-manifold.   By the classification theory of ALE
hyper-K\"ahler 4-manifold (c.f. \cite{Kr}),  there is a close
surface $\Sigma \subset N$ such that $[\Sigma]\in H_{2}(N,
\mathbb{Z})$, and $[\Sigma]\cdot[\Sigma]=-2$. Then, for $k\gg 1$ and
$r\gg 1$, $F_{k,r}(\Sigma)$ is a cycle in $M$, and defines  a
homology class $[F_{k,r}(\Sigma)]\in H_{2}(M, \mathbb{Z})$ with
$[F_{k,r}(\Sigma)]\cdot[F_{k,r}(\Sigma)]=-2$. Let $a\in \mathbb{Z}$
and  $b\in \mathbb{Z}$ such that $[F_{k,r}(\Sigma)]=aH+bE$, where
$H$ and $E$ are the  two generators of $H_{2}(M, \mathbb{Z})$ such
that $H\cdot H=1$, $E\cdot E=-1$ and $H\cdot E=0$. However,  the
equation $a^{2}-b^{2}=[F_{k,r}(\Sigma)]\cdot[F_{k,r}(\Sigma)]=-2$
does  not have integer solutions.  It is a contradiction. Thus $
\pi_{1}(N)\neq \{1\}$.
\end{proof}
è
 Actually, by the same arguments as  in the proof of (2) in Theorem
 5.2 of \cite{T},  we can see that  $
\pi_{1}(N)$ is a cyclic finite group, and $(N, J_{\infty},
g_{\infty})$ is asymptotic to $\mathbb{C}^{2}/\Gamma$, where
$\Gamma$ is a finite cyclic  subgroup of $U(2)$ given by Lemma 5.5
in \cite{T}. If one wants  to use the technique in the proof of
Theorem 1.1 to prove the bounding  of curvatures along the
K\"ahler-Ricci flow on $\mathbb{CP}^{2}\sharp
\overline{\mathbb{CP}}^{2}$, a method must be found to prove that
$N$ is actually simply  connected. In a recently preprint \cite{CW},
it is  claimed  that this  could be done by using the fact that
$\mathbb{CP}^{2}\sharp \overline{\mathbb{CP}}^{2}$ is  a toric
manifold. However, the full  details of the arguments in  \cite{CW}
  have not  appeared yet.

Our method can  be used to give a different proof of  the following
theorem in term of Gromov-Hausdorff convergence,
 which is already implied by \cite{TZ} where the Monge-Amp\`{e}re
 flow was used.

  \begin{theorem} \rm (\cite{TZ}) \it If $M$ is   holomorphic to
    $\mathbb{CP}^{2}$, and  $g(t)$, $t\in [0, +\infty)$, is   a solution of the normalized
K\"ahler-Ricci flow (\ref{1.1}),  then, for any  sequence of times
$t_{k}\rightarrow \infty$, a subsequence of $(M, g(t_{k}))$
converges smoothly to the unique K\"ahler-Einstein metric on
$\mathbb{CP}^{2}$ in the Cheeger-Gromov sense.
   \end{theorem}

   \begin{proof} It is well known  that,  on
$\mathbb{CP}^{2}$,  there is a unique K\"ahler-Einstein metric
presenting $2\pi c_{1}$, the Fubini-Study metric (c.f. \cite{BM}).
This implies that the Mabuchi's K-energy $\nu_{g_{0}}(g(t))$ is
bounded from below (c.f. \cite{PS}). Thus,  by Perelman's estimate
for scalar curvatures,  (6.1) in \cite{PS} holds, i.e.
$$\int_{M}|\partial\overline{\partial}u_{t}|^{2}dv_{t}\longrightarrow
0,$$ when $t\rightarrow \infty$, where $u_{t}$ are functions
satisfying
$-Ric(g(t))+g(t)=\sqrt{-1}\partial\overline{\partial}u_{t}$. Since
the Hodge Laplacian satisfies
$\triangle=2(\partial^{*}\partial+\partial\partial^{*})=2(
\overline{\partial}^{*}\overline{\partial}+\overline{\partial}\overline{\partial}^{*})$,
$\triangle \overline{\partial}=\overline{\partial} \triangle$, and
$\triangle u_{t}=R(g(t))-4$, we obtain that
\begin{equation}\int_{M}|R(g(t))-4|^{2}dv_{t}=\int_{M}|\triangle u_{t}|^{2}dv_{t}
=4\int_{M}|\partial\overline{\partial}u_{t}|^{2}dv_{t}\longrightarrow
0. \end{equation} Then
\begin{eqnarray*} \lim\limits_{t\rightarrow \infty}\int_{M}R(g(t))^{2}dv_{t} &
= &  \lim\limits_{t\rightarrow \infty}\int_{M}(8R(g(t))-16)dv_{t} \\
& = & \lim\limits_{t\rightarrow \infty}(16\int_{M}Ric(g(t))\wedge
g(t)- 8\int_{M}g(t)\wedge g(t))\\ & = &
32\pi^{2}\int_{M}c_{1}^{2}(M)=32\pi^{2}(2\chi(M)+3\tau(M)).
\end{eqnarray*} By (3.2), we have $$\lim\limits_{t\rightarrow
\infty}\int_{M}|W^{-}(g(t))|^{2}dv_{t}=\frac{8}{3}\pi^{2}(\chi(M)-3\tau(M))=0,$$
since $\chi(M)=3$ and $\tau(M)=1$.

 If   $\sup\limits_{t\in [0,
\infty)}|Rm(g(t))|=+\infty,$ then  there exists  a sequence of times
$t_{k}\rightarrow\infty$, and a sequence of  points $x_{k}\in M$
such that
$$Q_{k}=|Rm(g(t_{k}))|(x_{k})=\sup_{M\times[0,t_{k}]}|Rm(g(t))|\rightarrow\infty.$$
Let $\widetilde{g}_{k}=Q_{k}g(t_{k})$.  By  the same arguments as in
the proof of Theorem 1.1, by passing to a subsequence, $\{(M, J,
\widetilde{g}_{k}, x_{k})\}$ converges smoothly to a  complete
Ricci-flat K\"ahler surface   $(N, J_{\infty}, g_{\infty},
x_{\infty})$ with $$ \sup_{N} |Rm(g_{\infty})|=1.$$ Furthermore, by
the smooth convergence,
$$\int_{N}|W^{-}(g_{\infty})|^{2}dv_{\infty}\leq \lim\limits_{t\rightarrow
\infty}\int_{M}|W^{-}(g(t))|^{2}dv_{t}=0, \ \ {\rm thus} \ \
W^{-}(g_{\infty})\equiv 0,$$ on $N$. From
$R^{2}(g_{\infty})=24|W^{+}(g_{\infty})|^{2}\equiv 0$, we obtain
that $|Rm(g_{\infty})|\equiv 0$ on $N$. It is a contradiction. Hence
there is a constant $C>0$ independent of $t$  such that $$
|Rm(g(t))|\leq C.$$

Finally, by   the same arguments as in the proof of Theorem 1.1, for
any $t_{k}\rightarrow \infty$, a subsequence of $( M, g(t_{k}+t))$,
$t\in [0, 1]$, converges smoothly to a compact K\"ahler-Ricci
soliton $(X, h(t))$, $t\in [0, 1]$. By (3.3), $R(h(t))\equiv 4$,
and, thus, $h(0)$ is a K\"ahler-Einstein metric.  By Kodaria
classification theorem, it is well known that the Fano surface
diffeomorphic to $\mathbb{CP}^{2}$ is unique. Therefore, $X\cong M$.
   \end{proof}

\section{Proof of Theorem    1.2}

The proof of Theorem  1.2 follows by a similar argument as in the
proof of Theorem 1.1.

\begin{proof}[Proof of Theorem    1.2]
Suppose not, there exist a sequence of  $t_{k}\rightarrow T$ and
points $x_{k}\in M$ such that
$$Q_{k}=|Rm(g(t_{k}))|(x_{k})=\sup_{M\times[0,t_{k}]}|Rm(g(t))|\rightarrow\infty.$$
Then consider the sequence of solutions to the Ricci flow
$$(M,Q_{k}g(Q_{k}^{-1}t+t_{k}),x_{k}),t\in[-Q_{k}t_{k},0].$$
First we assume $g(t)$ is a solution to the unnormalized Ricci flow.
Perelman's no local collapsing theorem (cf. \cite[Theorem 4.1]{Pe}
or \cite[Remark 12.13]{KL}) applies to show that there is a
$\kappa>0$ such that $Vol(B(x,r,g(t)))\geq\kappa r^{n}$ for each
metric ball $B(x,r,g(t))$ in $(M,g(t))$ with radius $r\leq\sqrt{T}$.
Then using Hamilton's compactness theorem for Ricci flow solutions,
the sequence will converge modulo a subsequence to another solution
to the Ricci flow, say
$(M_{\infty},g_{\infty}(t),x_{\infty}),t\in(-\infty,0]$, which has
the properties that $Vol(B(r,g_{\infty}(t)))\geq\kappa r^{n},$ for
each metric ball $B(r,g_{\infty}(t))$ in
$(M_{\infty},g_{\infty}(t))$ of radius $r$, and that
$$\int_{M_{\infty}}|Rm(g_{\infty}(t))|^{n/2}dv_{g_{\infty}(t)}\leq\limsup_{k\rightarrow\infty}\int_{M}|Rm(g(t))|^{n/2}dv_{t}<\infty,$$
$$|Rm(g_{\infty}(0))|(x_{\infty})=1\mbox{ and }R(g_{\infty}(t))\equiv0\mbox{ over }M_{\infty}\times(-\infty,0].$$
From the evolution of the volume $Vol(g(t))$ of the metric $g(t)$:
$$\frac{d}{dt}Vol(g(t))=-\int_{M}R(g(t))dv_{g(t)}\geq-CVol(g(t)),$$
we conclude that $Vol(g(t))\geq Vol(g(0))e^{-Ct}\geq
Vol(g(0))e^{-CT}$ for each metric $g(t)$. So the limits
$(M_{\infty},g_{\infty}(t))$ are non-compact Ricci flat manifolds.
After a double covering, we may also assume that the manifold $M$ is
oriented and so the limit $M_{\infty}$ is also oriented. By odd
dimensional assumption of $M$, using theorem 3.5 of \cite{An1}, we
conclude that $(M_{\infty},g_{\infty}(0))$ is in fact the Euclidean
space, which contradicts the fact
$|Rm(g_{\infty}(0))|(x_{\infty})=1$.

If $g(t)$ is a solution to the normalized Ricci flow, then the
rescaling factor from $g(t)$ to the corresponding unnormalized Ricci
flow is uniformly bounded from above and below (stays bounded away
from zero), since the scalar curvature of $g(t)$ is absolutely
bounded. Thus the corresponding unnormalized Ricci flow exists in
finite time, and then Perelman's no local collapsing theorem uses
also. So repeatedly, $Vol(B(x,r,g(t)))\geq\kappa r^{n}$ for each
metric ball $B(x,r,g(t))$ in $(M,g(t))$ with radius $r\leq\sqrt{T}$,
for some universal $\kappa>0$. If $g(t)$ does not has uniformly
bounded Riemannian curvature, then there is a sequence of times
$t_{k}\rightarrow T$ and points $x_{k}\in M$ such that
$$Q_{k}=|Rm(g(t_{k}))|(x_{k})=\sup_{M\times[0,t_{k}]}|Rm(g(t))|\rightarrow\infty.$$
Then consider the sequence of solutions to the Ricci flow
$$(M,Q_{k}g(Q_{k}^{-1}t+t_{k}),x_{k}),t\in[-Q_{k}t_{k},0],$$
which will converge along a subsequence to a Ricci flat solution on
an open manifold. The limit solution is flat by a same argument and
we obtain a contradiction.
\end{proof}

\end{document}